\documentclass[12pt,reqno]{amsart}

\usepackage{amssymb}
\usepackage{amsfonts}

\newcommand{\nc}{\newcommand}
\nc{\thref}[1]{Theorem~\ref{theo:#1}}
\nc{\selabel}[1]{\label{sect:#1}}
\nc{\seref}[1]{Section~\ref{sect:#1}}
\nc{\lelabel}[1]{\label{lemm:#1}}
\nc{\leref}[1]{Lemma~\ref{lemm:#1}}
\nc{\prlabel}[1]{\label{prop:#1}}
\nc{\prref}[1]{Proposition~\ref{prop:#1}}
\nc{\colabel}[1]{\label{coro:#1}}
\nc{\coref}[1]{Corollary~\ref{coro:#1}}
\nc{\exlabel}[1]{\label{exam:#1}}
\nc{\exref}[1]{Example~\ref{exam:#1}}
\nc{\delabel}[1]{\label{defi:#1}}
\nc{\deref}[1]{Definition~\ref{defi:#1}}
\nc{\eqlabel}[1]{\label{equa:#1}}
\nc{\relabel}[1]{\label{rema:#1}}
\nc{\reref}[1]{Lemma~\ref{rema:#1}}
\providecommand{\operatorname}[1]{\mathrm{#1}\,}
\nc{\Hom}{\operatorname{Hom}}
\nc{\Mor}{\operatorname{Mor}}
\nc{\Aut}{\operatorname{Aut}}
\nc{\Ann}{\operatorname{Ann}}
\nc{\Ker}{\operatorname{Ker}}
\nc{\Trace}{\operatorname{Trace}}
\nc{\Char}{\operatorname{Char}}
\nc{\Mod}{\operatorname{Mod}}
\nc{\End}{\operatorname{End}}
\nc{\Spec}{\operatorname{Spec}}
\nc{\Span}{\operatorname{Span}}
\nc{\sgn}{\operatorname{sgn}}
\nc{\Id}{\operatorname{Id}}
\nc{\Com}{\operatorname{Com}}
\nc{\rank}{\operatorname{rank}}

\let\dirsum=\oplus

\let\:=\colon

\newtheorem{de}{Definition}[section]
\newtheorem{lm}[de]{Lemma}
\newtheorem{pr}[de]{Proposition}
\newtheorem{co}[de]{Corollary}
\newtheorem{re}[de]{Remark}
\newtheorem{res}[de]{Remarks}
\newtheorem{te}[de]{Theorem}
\newtheorem{ex}[de]{Example}
\newtheorem{exs}[de]{Examples}

\def\bex{\begin{ex}}
\def\eex{\end{ex}}
\def\bexs{\begin{exs}}
\def\eexs{\end{exs}}
\def\bl{\begin{lm}}
\def\el{\end{lm}}
\def\bc{\begin{co}}
\def\ec{\end{co}}
\def\bt{\begin{te}}
\def\et{\end{te}}
\def\bpr{\begin{pr}}
\def\epr{\end{pr}}
\def\br{\begin{re}}
\def\er{\end{re}}
\def\brs{\begin{res}}
\def\ers{\end{res}}
\def\bd{\begin{de}}
\def\ed{\end{de}}
\def\be{\begin{equation}}
\def\ee{\end{equation}}
\def\bea{\begin{eqnarray*}}
\def\eea{\end{eqnarray*}}
\def\bp{\begin{proof}}
\def\ep{\end{proof}}

\def\l{{\lambda}}
\def\L{{\Lambda}}
\def\g{{\gamma}}
\def\G{{\Gamma}}
\def\d{{\delta}}
\def\D{{\Delta}}
\def\s{{\sigma}}

\def\ZZ{{\mathbb Z}}
\def\NN{{\mathbb N}}

\begin{document}
\title{Valuation Extensions of Filtered and Graded Algebras}
\author{C. Baetica$^1$$^,$$^*$ and F. Van Oystaeyen$^2$\\ \\
{$^1$University of Bucharest, Faculty of Mathematics, Bucharest,
Romania}\\ \\ {$^2$University of Antwerp, Department of
Mathematics and Computer Science, Antwerp, Belgium}}

\begin{abstract}
In this note we relate the valuations of the algebras appearing in
the non-commutative geometry of quantized algebras to properties
of sub-lattices in some vector spaces. We consider the case of
algebras with $PBW$-bases and prove that under some mild
assumptions the valuations of the ground field extend to a
non-commutative valuation. Later we introduce the notion of
$F$-reductor and graded reductor and reduce the problem of finding
an extending non-commutative valuation to finding a reductor in an
associated graded ring having a domain for its reduction.\\



\end{abstract}

\thanks{$^*$Correspondence: C. Baetica, University of Bucharest,
Faculty of Mathematics, Str. Academiei 14, RO--010014, Bucharest,
Romania; E--mail: baetica@al.math.unibuc.ro}
\date{}
\maketitle

\pagestyle{myheadings}

\markboth{\sc C. Baetica and F. Van Oystaeyen} {\sc Valuation
Extensions}
\section*{Introduction}
\noindent In the so-called non-commutative geometry of quantized
algebras, some steps in the direction of a divisor theory via
non-commutative valuation rings have already been taken; see
Hussein and Van Oystaeyen (1996), Van Oystaeyen and Willaert
(1996), Willaert (1997), etc. Most of the algebras appearing in
the theory of rings of differential operators, quantized algebras
of different kinds including many quantum groups, regular algebras
in projective non-commutative geometry, come equipped with a
natural gradation or filtration controlled by some finite
dimensional vector space(s), e.g. the degree one part of
filtration or gradation.

In this note we relate the valuations of the algebras considered
to properties of sub-lattices in some vector space(s) as mentioned
above. In the first section we generalize some results from Li and
Van Oystaeyen (1990) to $\G$-filtrations, where $\G$ is a totally
ordered group. In Section 2 we consider the case of algebras with
$PBW$-bases and prove that under some mild assumptions the
valuations of the ground field extend to a non-commutative
valuation. As an application we single out the case of Weyl
algebras, although the results of this section apply to a large
class of examples. In the last two sections we introduce the
notion of $F$-reductor and graded reductor and reduce the problem
of finding an extending non-commutative valuation to finding a
reductor in an associated graded ring having a domain for its
reduction.

Throughout this paper we use the notation in Van Oystaeyen (2000).

\section{Valuation filtrations}
\noindent Let $R$ be a ring and $\G$ a totally ordered group. A
family $FR$ of additive subgroups $F_{\g}R$, $\g\in\G$ satisfying
\begin{itemize}
\item[$(i)$] $\g\leq\d$ implies $F_{\g}\subseteq F_{\d}$,
\item[$(ii)$] $F_{\g}R\;F_{\d}R\subseteq F_{\g+\d}R$, for
all $\g,\d\in\G$,
\item[$(iii)$] $1\in F_0R$,
\end{itemize}
is a {\em filtration of type $\G$}, or a {\em $\G$-filtration} on
$R$. For a $\G$-filtration $FR$ we may define
the associated graded ring
$G_F(R)=\dirsum_{\g\in\G}F_{\g}R/F_{<\g}R$, where
$F_{<\g}R=\sum_{\g'<\g}F_{\g'}R$.

We say that $FR$ is {\em $\G$-separated} if for every $a\in R$,
$a\neq 0$ there is a $\g\in\G$ such that $a\in F_{\g}R-F_{<\g}R$.
For $\G=\ZZ$ the $\G$-separatedness is equivalent to
$\cap_{n\in\ZZ}F_nR=0$, see Van Oystaeyen (2000, Remark 4.2.3),
but for an arbitrary $\G$ the latter condition may be strictly
weaker than this. If $FR$ is $\G$-separated, then we may define
the {\em principal symbol map} $\sigma:R\longrightarrow G_F(R)$ by
$\sigma(a)=a$ mod $F_{<\g}R$ whenever $a\in F_{\g}R-F_{<\g}R$. The
{\em degree} $\deg\sigma(a)=\g$ of $\sigma(a)$ is uniquely
determined. To $FR$ one associates a value-function
$v_F:R\longrightarrow \G\cup\{\infty\}$ defined by $v_F(0)=\infty$
and $v_F(a)=-\deg\sigma(a)$ for $a\neq 0$. It is well known that $v_F$
is a valuation function on $R$ whenever $G_F(R)$ is a domain; see
Van Oystaeyen (2000, Corollary 4.2.7). Furthermore, we have that
$1\in F_0R-F_{<0}R$, or equivalently $G_F(R)_0\neq 0$.

All the filtrations we are dealing with are considered to be {\em
exhaustive}, that is $\cup_{\g\in\G}F_{\g}R=R$.

\bd If $\D$ is a skewfield and $\L\subset\D$ is a subring, then
$\L$ is a {\em valuation ring} of $\D$ if it is invariant under
inner automorphisms of $\D$ and for $x\in\D$, $x\neq 0$, either
$x\in\L$ or $x^{-1}\in\L$.\ed

\bpr\label{characterization} Let $R$ be an Artinian ring and $FR$
a $\G$-separated filtration on $R$. Then the following are
equivalent:
\begin{itemize}
\item[$(i)$] $G_F(R)$ is a domain,
\item[$(ii)$] $R$ is a skewfield and $G_F(R)$ is a gr-skewfield,
\item[$(iii)$] $R$ is a skewfield, $F_0R$ is a valuation ring of
$R$ (with $F_{<0}R$ as unique maximal ideal), $\G_s=\{\g\in\G:
G_F(R)_{\g}\neq 0\}$ is a subgroup of $\G$ and $G_F(R)$ is
$\G_s$-strongly graded.
\end{itemize} \epr \bp (i) $\Rightarrow$ (ii) If $G_F(R)$ is a domain,
then $R$ is a domain; see Van Oystaeyen (2000, Proposition 4.2.4).
Therefore $R$ is a skewfield, because it is Artinian. For an
element $\sigma(a_{\g})\in G_F(R)_{\g}$, we have $a_{\g}\in
F_{\g}R-F_{<\g}R$ and hence $a_{\g}$ is invertible. Stand $a_{\d}$
for its inverse. From the assumption $1\in F_0R-F_{<0}R$, it
follows that $\g+\d\geq 0$. Necessarily we get $\g+\d=0$,
otherwise $\sigma(a_{\g})\sigma(a_{\d})=0$ contradicts the
assumption that $G_F(R)$ is a domain. Thus $\d=-\g$ and
$\sigma(a_{\d})$ is an inverse of $\sigma(a_{\g})$.\\ (ii)
$\Rightarrow$ (iii) Since $G_F(R)$ is a gr-skewfield, it is in
particular a gr-domain. As $\G$ is a totally ordered group,
$G_F(R)$ is a domain. By using again Van Oystaeyen (2000,
Proposition 4.2.4) we deduce that the value-function $v_F$
associated to the filtration $FR$ is a valuation function on $R$
whose valuation ring is $F_0R$. The fact that $\G_s$ is a subgroup
of $\G$ and $G_F(R)$ is $\G_s$-strongly graded is obvious. \\
(iii) $\Rightarrow$ (i) $G_F(R)$ is strongly graded and
$G_F(R)_0=F_0R/F_{<0}R$ is a skewfield, so $G_F(R)$ is in
particular a domain. \ep

\br {\em If $FR$ is a $\G$-separated {\em strong filtration}, i.e.
$F_{\g}R\;F_{\d}R=F_{\g+\d}R$ for all $\g,\d\in\G$, then
$\G_s=\G$. In particular, the valuation function $v_F$ is
surjective.}\er

By Proposition \ref{characterization} we get the following

\bc If $\Delta$ is a skewfield and $\L\subset\Delta$ is a subring,
then the following are equivalent:
\begin{itemize}
\item[$(i)$] $\L$ is a valuation ring of $\Delta$,
\item[$(ii)$] There is a totally ordered group $\G$ and
a valuation $v:\Delta\longrightarrow\G\cup\{\infty\}$ such that
$\L$ is the valuation ring of $v$,
\item[$(iii)$] There is a totally ordered group $\G$ and $F\D$
an exhaustive $\G$-separated filtration on $\D$ such that
$\L=F_0\D$ and $G_F(\D)$ is a domain.
\end{itemize} \ec \bp (i) $\Rightarrow$ (ii) See Van Geel (1981).\\
(ii) $\Rightarrow$ (iii) Set $F_{\g}\D=\{x\in\D: v(x)\geq-\g\}$,
for $\g\in\G$. Since $v$ is surjective, the filtration $F\D$ is
$\G$-separated. Furthermore, $F\D$ is a strong filtration and
thus $G_F(\D)$ is strongly graded with $G_F(\D)_0$ a skewfield.
This is enough in order to show that $G_F(\D)$ is a domain. \\
(iii) $\Rightarrow$ (i) Easily enough. \ep

\bc Let $FR$ be a $\G$-separated filtration on a ring $R$. Then
the following are equivalent:
\begin{itemize}
\item[$(i)$] $G_F(R)$ is a gr-skewfield and $F_{<0}R\subseteq J(F_0R)$,
\item[$(ii)$] $R$ is a skewfield, $F_0R$ is a valuation ring of
$R$ (with $F_{<0}R$ as unique maximal ideal), $\G_s=\{\g\in\G:
G_F(R)_{\g}\neq 0\}$ is a subgroup of $\G$ and $G_F(R)$ is
$\G_s$-strongly graded.
\end{itemize} \ec \bp (i) $\Rightarrow$ (ii) Let $a$ be a non-zero
element of $R$. Then there is $\g\in\G$ such that $\sigma(a)\in
G_F(R)_{\g}$. By hypothesis $G_F(R)$ is a gr-skewfield, hence
there is $b\in R$ with $\sigma(a)\sigma(b)=1$, that is $ab-1\in
F_{<0}R$. It means that $ab$ is invertible, hence $a$ is also
invertible. The other claims follow by Proposition
\ref{characterization}.\\ (ii) $\Rightarrow$ (i) Readily by
Proposition \ref{characterization}. \ep

\section{$K$-algebras with PBW-bases}
\noindent From the preceding section we have seen how important is
the $\G$-separatedness for $\G$-filtrations. A wide class of
algebras are known to possess $\G$-separated filtrations and
among them we mention the $K$-algebras with a finite $PBW$-basis.

Throughout this section we consider $K$ as being a field, $\G$ a
totally ordered group and $v:K\longrightarrow \G\cup\{\infty\}$ a
valuation on $K$. We always assume that $v$ is surjective, hence
$\G$ is a commutative group. Set $O_v$ the valuation ring of $K$
associated to $v$, $m_v$ its unique maximal ideal and
$k_v=O_v/m_v$ the residue field. Then we define a $\G$-filtration
$f^vK$ on $K$ by $f^v_{\g}K=\{x\in K: v(x)\geq -\g\}$ and we call
it the {\em valuation filtration} on $K$. This is an exhaustive
$\G$-separated strong filtration.

Let $A$ be an affine $K$-algebra generated by $a_1,\ldots,a_n$,
$K<\underline{X}>$ the free $K$-algebra on the set
$\underline{X}=\{X_1,\ldots,X_n\}$ and
$\pi:K<\underline{X}>\longrightarrow A$ the canonical $K$-algebras
morphism given by $\pi(X_i)=a_i$, $i=1,\ldots,n$. Restriction of
$\pi$ to $O_v<\underline{X}>$ defines a subring $\L$ of $A$, i.e.
$\L=\pi(O_v<\underline{X}>)$. Denote by $\mathcal{R}$ the ideal of
relations of $A$ which is obviously a two-sided ideal. If
$\mathcal{R}$ is generated by
$p_1(\underline{X}),\ldots,p_d(\underline{X})$ as a left ideal,
then it is generated by
$p_1(\underline{X}),\ldots,p_d(\underline{X})$ as a two-sided
ideal; see Hussein and Van Oystaeyen (1996, Lemma 2.1). Then we
may assume that each $p_i(\underline{X})\in O_v<\underline{X}>$
(up to multiplying by a suitable constant), but not all are in
$m_v<\underline{X}>$. Since $\mathcal{R}\cap O_v<\underline{X}>$
need not be generated again by the
$p_1(\underline{X}),\ldots,p_d(\underline{X})$ as a two-sided
ideal of $O_v<\underline{X}>$, we may encounter reduction problems.
We say that $\L$ defines a {\em good reduction} of $A$ (or that
$A$ {\em reduces well} at $O_v$) whenever $\overline{\mathcal{R}}$
is generated (as a two-sided ideal of $k_v<\underline{X}>$) by the
residues $\overline{p_i(\underline{X})}$, $i=1,\ldots,d$. For
instance, if
$\L'=O_v<\underline{X}>/<p_1(\underline{X}),\ldots,p_d(\underline{X})>$
is a flat (equivalently, torsion-free) $O_v$-module, then $\L$ is
a good reduction of $A$.

The subring $\L$ yields a filtration $F^vA$ on $A$ given by
$F^v_{\g}A=(f^v_{\g}K)\L$, $\g\in\G$. We call $F^vA$ the {\em
valuation filtration} on $A$. Obviously, the valuation filtration
$F^vA$ is exhaustive. By $G_v(A)$ we denote the associated graded
ring to the valuation filtration $F^vA$.

The next result is Theorem 3.4.7 from Van Oystaeyen (2000).

\bl\label{PBW} For a graded $K$-algebra $A$ that has a finite
$PBW$-basis, the valuation filtration $F^vA$ is $\G$-separated and
strong.\el

This enables us to extend $v$ to every graded $K$-algebra $A$ that
has a finite $PBW$-basis. In the following we denote by $A(k)$ the
affine $k$-algebra defined by the same relations as $A$, where $k$
is a suitable field.

\bt\label{extension-PBW} Let $A$ be a graded $K$-algebra with a
finite $PBW$-basis and $\L\subset A$ as before.
\begin{itemize}
\item[$(i)$] If $\overline{\L}=\L/m_v\L$ is a domain and $A$ is an
Ore domain, then every valuation $v$ on $K$ extend to
$Q=Q_{cl}(A)$, the classical ring of fractions of $A$.
\item[$(ii)$] If $\L$ defines a good reduction of $A$, $A(k_v)$
is a domain and $A$ is an Ore domain, then every valuation $v$ on
$K$ extend to $Q$.
\end{itemize} \et \bp (i) The associated graded ring $G_v(A)$ is
strongly graded with $G_v(A)_0=\overline{\L}$, hence $G_v(A)$ is a
domain. By using Van Oystaeyen (2000, Corollary 4.2.7) we get a
valuation function on $A$. It turns out that this is an extension
of $v$ to $A$ and then we extend the valuation function on $A$ to
$Q$ in the usual manner. All we have to do now is to show that
$F^v_{\g}A\cap K=f^v_{\g}K$. It is clear that $f^v_{\g}K\subset
F^v_{\g}A\cap K$. For converse, pick an element $x\in
F^v_{\g}A\cap K$. Thus $x\in(f^v_{\g}K)\L$ and $x\in K$, therefore
$(f^v_{-\g}K)x\subset\L\cap K$. As $\L\cap K=O_v$, we get
$(f^v_{\g}K)(f^v_{-\g}K)x\subset f^v_{\g}K$. But
$(f^v_{\g}K)(f^v_{-\g}K)=f^v_0K=O_v$ and it is enough to conclude
that $x\in f^v_{\g}K$. \\ (ii) Since $\L$ defines a good reduction
of $A$, the ring $\overline{\L}$ is isomorphic to $A(k_v)$, hence
it is a domain. Now we can apply (i). \ep

\bc Let $D_n(K)$ be the skewfield of fractions of the Weyl algebra
$A_n(K)$. Then every valuation $v$ on $K$ extends to $D_n(K)$. \ec
\bp Set $A=A_n(K)$. Then $A$ satisfies the conditions of Lemma
\ref{PBW} and thus the valuation filtration $F^vA$ is
$\G$-separated. The associated graded ring $G_v(A)$ is then
strongly graded with $G_v(A)_0=\L/m_v\L$, where $\L$ is defined as
above. Note that $\L$ has also a finite $PBW$-basis over $O_v$. In
particular $\L$ is $O_v$-flat, therefore it defines a good
reduction of $A$. Consequently we can apply Proposition
\ref{extension-PBW}(ii) in order to get that every valuation $v$
on $K$ extend to $D_n(K)$. \ep

For graded $K$-algebras with a finite $PBW$-basis the associated
graded ring defined by the valuation filtration is a twisted group
ring.

\bpr\label{crossed} If $A$ is a graded $K$-algebra that has a
finite $PBW$-basis, then the associated graded ring $G_v(A)$ is
isomorphic to the twisted group ring $\overline{\L}*\G$, where
$\overline{\L}=\L/m_v\L$. \epr \bp The Rees algebra
$\widetilde{A}$ of $A$ associated to the valuation filtration
$F^vA$ is defined by $\widetilde{A}=\dirsum_{\g\in\G}F^v_{\g}A$.
Since $\L$ is $O_v$-flat, we obtain
$$\widetilde{A}=\dirsum_{\g\in\G}(f_{\g}K)\L\simeq
(\dirsum_{\g\in\G}f_{\g}K)\otimes_{O_v}\L=\widetilde{K}\otimes_{O_v}\L.$$
From $G_v(A)\simeq\widetilde{A}/\widetilde{A}\G^+$ we get that
$G_v(A)\simeq \widetilde{K}/\widetilde{K}\G^+\otimes_{O_v}\L\simeq
G_v(K)\otimes_{O_v}\L$. On the other side, $G_v(K)$ is a
commutative strongly graded ring with the trivial action,
therefore it is isomorphic to the twisted group ring $k_v*\G$ and
thus we get
$G_v(A)\simeq(k_v*\G)\otimes_{O_v}\L\simeq\L/m_v\L*\G$.\ep

In particular, for Weyl algebras we obtain

\bc $G_v(A_n(K))\simeq A_n(k_v)*\G$.\ec

For a plenty of examples of algebras with a finite $PBW$-basis,
the reader is referred to Berger (1992) and Van Oystaeyen (2000,
Observation 4.3.8).

Similarly to Hussein and Van Oystaeyen (1996, Proposition 2.3) we
can prove the following

\bpr Let $A$ be a graded $K$-algebra with a finite $PBW$-basis and
$\L\subset A$ as before. If $A$ is a prime Goldie ring and
$\overline{\L}$ is a domain, then
\begin{itemize}
\item[$(i)$] $Q=Q_{cl}(A)$ is a skewfield and $F^vA$ extends to a
strong filtration on $Q$ such that $F^v_0Q$ is a valuation ring of
$Q$ containing $O_v$, that is $v$ extends to $Q$.
\item[$(ii)$] $Q'=Q^g(A)$, the graded ring of fractions of $A$, is
a gr-skewfield and $F^vQ$ induces a graded valuation on $Q'$.
\end{itemize}
\epr

\section{Filtered and graded reductors}
\noindent Let $A$ be a $K$-algebra, $K$ a field and $FA$ an
exhaustive separated $\ZZ$-filtration such that $K\subseteq F_0A$.
Suppose in addition that $FA$ is {\em finite}, i.e.
$\dim_KF_nA<\infty$ for all $n\in\ZZ$. Since $FA$ is finite and
separated, it must be {\em left limited}, i.e. there is
$n_0\in\ZZ$ such that $F_nA=0$ for all $n\leq n_0$. Without loss
of generality we may suppose that the filtration $FA$ is {\em
positive}, that is $F_nA=0$ for all $n<0$.

Let $v$ be a valuation on $K$ with the value group $\G$,
$O_v\subset K$ the corresponding valuation ring, $m_v$ its unique
maximal ideal and $k_v=O_v/m_v$ the residue field.

\bd Let $V$ be a finite dimensional $K$-vector space and $M$ an
$O_v$-submodule of $V$. Then $M$ is an {\em $O_v$-lattice} in $V$
if it contains a $K$-basis of $V$ and is a submodule of a finitely
generated $O_v$-submodule of $V$. Usually, we will denote the
$k_v$-vector space $M/m_vM$ by $\overline{V}$. \ed

\br\label{rem} {\em For any $O_v$-lattice $M$ in $V$ we have
$\dim_{k_v}\overline{V}\leq\dim_KV$. When equality holds we say
that $M$ defines an {\em unramified reduction} of $V$.} \er

The next result deals with some elementary properties of
$O_v$-lattices that we did not find explicitly in the literature.

\bl\label{lattices} Let $V$ be a finite dimensional $K$-vector
space, $V'\subset V$ a $K$-vector subspace and $M\subset V$ an
$O_v$-submodule.
\begin{itemize}
\item[$(i)$] If $M$ is an $O_v$-lattice in $V$, then the quotient
module $M/M\cap V'$ is an $O_v$-lattice in $V/V'$.
\item[$(ii)$] If $M/M\cap V'$ is an $O_v$-lattice in $V/V'$ and
$M\cap V'$ is an $O_v$-lattice in $V'$, then $M$ is an
$O_v$-lattice in $V$.
\end{itemize} \el \bp (i) It is known that
$M\cap V'$ is an $O_v$-lattice in $V'$; see Fossum (1973,
Proposition 2.2(ii)). Then it is straightforward that $M/M\cap V'$
is an $O_v$-lattice in $V/V'$.\\ (ii) Easy. \ep

The {\em valuation filtration} $F^vV$ on $V$ defined by
$F^v_{\g}V=(f^v_{\g}K)M$ is an exhaustive filtration, but {\em we
do not know if it is $\G$-separated or not.} It is rather easy to
see that the valuation filtration is $\G$-separated whenever $M$
is a free $O_v$-module. This happens, for instance, when $\G=\ZZ$.
This shows that the discrete case can be handled with ease.

\bd Let $A$ be a filtered $K$-algebra with a finite filtration
$FA$ and $\L\subset A$ a subring. Then $\L$ is an {\em
$F$-reductor} if $\L\cap K=O_v$ and $\L\cap F_nA$ is an
$O_v$-lattice in $F_nA$ for all $n\in\NN$. \ed We call the ring
$\overline{A}=\L/m_v\L$ the (filtered) {\em reduction} of $A$ with
respect to $\L$. The {\em valuation filtration} $F^vA$ on $A$ is
defined by $F^v_{\g}A=(f^v_{\g}K)\L$, $\g\in\G$, while the {\em
induced filtration} $F\L$ of $FA$ on $\L$ is given by
$F_n\L=\L\cap F_nA$, $n\in\NN$.

\bpr\label{F-reductor} Let $A$ be as before and $\L\subset A$ an
$F$-reductor. We have
\begin{itemize}
\item[$(i)$] $A=K\L$ and the valuation filtration $F^vA$ is
exhaustive,
\item[$(ii)$] $F^vA$ is $\G$-separated, provided that
$F^v(F_nA)$ is $\G$-separated for all $n\in\NN$,
\item[$(iii)$] the induced filtration of $F^vA$ on $F_nA$ is the
valuation filtration on $F_nA$ and provided that it is $\G$-separated
we get $G_v(F_nA)=\overline{F_nA}*\G$, where
$\overline{F_nA}=F_n\L/m_vF_n\L$,
\item[$(iv)$] $(\overline{F_nA})_{n\in\NN}$ defines a finite filtration
on $\overline{A}$, denoted by $\overline{F}\;\overline{A}$ and
$G_{\overline{F}}(\overline{A})=\overline{G(A)}$, where
$\overline{G(A)}=G(\L)/m_vG(\L)$.
\end{itemize}\epr \bp (i) Since $F_n\L$ is an $O_v$-lattice in $F_nA$,
we have $K(F_n\L)=F_nA$ for all $n\in\NN$. As $FA$ is exhaustive
we get $K\L=A$. The valuation filtration $f^vK$ is exhaustive and
$K\L=A$, therefore $F^vA$ is an exhaustive filtration.\\ (ii) In
order to show that $F^vA$ is $\G$-separated we first note that
\begin{equation}\label{1}
(f^v_{\g}K)\L\cap F_nA=(f^v_{\g}K)(\L\cap F_nA)
\end{equation}
for all $\g\in\G$ and $n\in\NN$. It is obvious that
$(f^v_{\g}K)(\L\cap F_nA)\subseteq (f^v_{\g}K)\L\cap F_nA$. The
converse inclusion follows readily if one observes that the
filtration $f^vK$ is strong. To enter the detail, pick an element
$x\in (f^v_{\g}K)\L\cap F_nA$. It follows that
$(f^v_{-\g}K)x\subseteq \L\cap F_nA$ and therefore $x\in
(f^v_{\g}K)(\L\cap F_nA)$.

For any element $x\in A$, $x\neq 0$, there exists an $n\in\NN$
such that $x\in F_nA$. Since the valuation filtration on $F_nA$ is
assumed to be $\G$-separated, there exists an element $\g\in\G$
such that $x\in F^v_{\g}(F_nA)-F^v_{<\g}(F_nA)$. By using
(\ref{1}) we obtain that $x\in F^v_{\g}A-F^v_{<\g}A$. \\ (iii) By
using (\ref{1}) again we get that the induced filtration of $F^vA$
on $F_nA$ is the valuation filtration on $F_nA$. The crossed
product structure of $G_v(F_nA)$ can be seen in the same way as in
the proof of Proposition \ref{crossed}.\\ (iv) Let us remark that
$m_v(\L\cap F_nA)=m_v\L\cap F_nA$. This can also be obtained by
using (\ref{1}). Then $(\overline{F_nA})_{n\in\NN}$ defines a
finite filtration on $\overline{A}$; see Remark \ref{rem}. It is
easily seen that $G_{\overline{F}}(\overline{A})=\overline{G(A)}$.
\ep

It is a common strategy to deduce properties of $A$ from
properties of $G_F(A)$ whenever possible, so let us focus on the
graded situation for a moment.

\bd Let $R$ be an $\NN$-graded $K$-algebra with $K\subseteq R_0$.
Suppose that $R$ is {\em locally finite}, i.e. $\dim_KR_n<\infty$
for all $n\in\NN$ and let $\L\subset R$ be a graded subring. Then
$\L$ is called a {\em graded reductor} if $\L\cap K=O_v$ and
$\L\cap R_n$ is an $O_v$-lattice in $R_n$ for all $n\in\NN$. \ed

The ring $\overline{R}=\L/m_v\L$ is called the (graded) {\em
reduction} of $R$ with respect to $\L$. The valuation filtration
$F^vR$ on $R$ is similarly defined by $F^v_{\g}R=(f^v_{\g}K)\L$,
$\g\in\G$. A $\G$-filtration $FR$ on $R$ such that
$F_{\g}R=\dirsum_{n\in\NN}F_{\g}R\cap R_p$, for all $\g\in\G$ is
called a {\em graded filtration}.


\bpr\label{gr-reductor} Let $R$ be as before and $\L\subset R$ a
graded reductor. We have
\begin{itemize}
\item[$(i)$] $R=K\L$ and the valuation filtration $F^vR$ is
an exhaustive graded filtration,
\item[$(ii)$] $F^vR$ is $\G$-separated, provided that $F^vR_n$
is $\G$-separated for all $n\in\NN$,
\item[$(iii)$] the induced filtration of $F^vR$ on $R_n$ is the
valuation filtration on $R_n$ and provided that it is
$\G$-separated we get $G_v(R_n)=\overline{R}_n*\G$, where
$\overline{R}_n=\L_n/m_v\L_n$,
\item[$(iv)$] $\overline{R}=\dirsum_n\overline{R}_n$ and
$G_v(R)=\overline{R}*\G$.
\end{itemize}\epr
\bp For (i) and (iii) we can argue similarly to the proof of
Proposition \ref{F-reductor}, while (iv) is rather obvious. \\
(ii) Pick an element $x\in R$, $x\neq 0$ and write
$x=\sum_{i=0}^nx_i$ with $x_i\in R_i$ not all zero. For each $x_i$
with $x_i\neq 0$ we get an element $\g_i\in\G$ such that $x_i\in
F_{\g_i}^vR_i-F_{<\g_i}^vR_i$. Assume that $\g$ is the maximum of
the $\g_i$. We claim that $x\in F_{\g}^vR-F_{<\g}^vR$. If $x\in
F_{<\g}^vR$, then there is $\d\in\G$, $\d<\g$, with the property
that $x\in F_{\d}^vR$. Since $F_{\d}^vR=(f_{\d}^vK)\L$, we can
write $x=\sum_{j=1}^s\alpha_j\l_j$, where $\alpha_j\in f_{\d}^vK$
and $\l_j\in\L$. Each $\l_j$ can be written as follows
$\l_j=\sum_n\l_{jn}$, with $\l_{jn}\in\L_n$, and thus we get
$x=\sum_n(\sum_{j=1}^s\alpha_j\l_{jn})$. In particular,
$x_i=\sum_{j=1}^s\alpha_j\l_{ji}$ and therefore
$x_i\in(f_{\d}^vK)\L_i$, a contradiction. \ep

Note that the reductor property for $\L$ is "almost" sufficient to
make $F^vR$ into a $\G$-separated filtration, whereas in Hussein
and Van Oystaeyen (1996) this was obtained (for $\G=\ZZ$) by
restricting to {\em connected} positively graded $K$-algebras,
that is, $R=K\dirsum R_1\dirsum\cdots\dirsum
R_n\dirsum\cdots=K[R_1]$ and $\dim_KR_1<\infty$.

\bl If $\L$ is a graded reductor for $R$ and $S$ is an Ore set in
$R$, then there is an Ore set $S'$ in $\L$ such that
$(S')^{-1}\L=S^{-1}R$. \el \bp Take $S'=K^{\times}S\cap\L$, where
$K^{\times}=K-\{0\}$. Since $K\L=R$, for every finite subset
$\{a_1,\ldots,a_m\}$ of $R$ there exists a nonzero $\lambda\in
O_v$ such that $\lambda a_i\in\L$ for all $i=1,\ldots,m$. It is
straightforward to check the Ore condition for $S'$ and moreover
$(S')^{-1}\L=S^{-1}R$ is clear as well. \ep

\bc Suppose that $\G=\ZZ$ and $G_v(R)$ is a domain. Then the
localized filtration on $S^{-1}R$ derived from $F^vR$, denoted by
$F^vS^{-1}R$, is exactly the localized filtration deriving from
$m_v\L$-adic filtration on $\L$. \ec

Similarly to Hussein and Van Oystaeyen (1996, Proposition 2.3) we
have the following

\bpr\label{ring of fractions} Let $R$ be as before and $\L\subset
R$ a graded reductor such that $F^vR$ is $\G$-separated. Suppose
that $R$ is a prime Goldie ring and $\overline{R}$ is a domain.
Then
\begin{itemize}
\item[$(i)$] $Q=Q_{cl}(R)$ is a skewfield and $F^vR$ extends to a
strong filtration on $Q$ such that $F^v_0Q$ is a valuation ring of
$Q$ containing $O_v$, that is $v$ extends to $Q$.
\item[$(ii)$] $Q'=Q^g(R)$, the graded ring of fractions of $R$, is
a gr-skewfield and $F^vQ$ induces a graded valuation on $Q'$.
\end{itemize}
\epr

For connected positively graded $K$-algebras, as are all the
graded algebras appearing in the non-commutative geometry of
$Proj$, the existence of graded reductors can be expressed in
terms of certain unramifiedness properties.

Let $R$ be a connected positively graded $K$-algebra with $\dim
R_1=n$, $\underline{X}=\{X_1,\ldots,X_n\}$, and take
$\pi:K<\underline{X}>\longrightarrow R$ a presentation of $R$. If
we set $\L=\pi(O_v<\underline{X}>)$, then
$\dim_{k_v}\overline{R}_1=n$ since no elements of degree one in
the gradation of $K<\underline{X}>$ are in $\mathcal{R}$, the
ideal of relations of $R$. Nevertheless, $\dim_KR_n$ and
$\dim_{k_v}\overline{R}_n$ may be different for $n>1$.

\bpr Let $R$ be a connected positively graded $K$-algebra.
\begin{itemize}
\item[$(i)$] If $\dim_KR_n=\dim_{k_v}\overline{R}_n$ for all $n\in\NN$,
then $\L$ is a graded reductor of $R$.
\item[$(ii)$] For $\G=\ZZ$ the converse of (i) holds.
\item[$(iii)$] When moreover $\overline{R}$ is a Goldie domain,
we have that $R$ is a domain and the $m_v\L$-adic filtration on
$R$ is induced by a valuation filtration on the skewfield of
microfractions of $R$.
\end{itemize}
\epr \bp (i) We have to prove that $\L_n\subset R_n$ is an
$O_v$-lattice for all $n\in\NN$. Take
$\{\overline{x}_1,\ldots,\overline{x}_m\}$ a $k_v$-basis of
$\overline{R}_n$. Then the elements $x_1,\ldots,x_m$ are linearly
independent over $K$ and by hypothesis they form a $K$-basis of
$R_n$. This shows that $\L_n$ contains a $K$-basis of $R_n$. As
$\L_n$ is a finitely generated $O_v$-module we can conclude that
$\L_n$ is an $O_v$-lattice in $R_n$.\\ (ii) When $\G=\ZZ$, each
graded component $\L_n$ of $\L$ is a free $O_v$-module whose rank
equals $\dim_KR_n$, as long as $\L$ is a graded reductor of $R$.\\
(iii) In this case $m_v=\pi O_v$, $\pi\in O_v$ and
$F_n^vR=\pi^{-n}\L$ for all $n\in\ZZ$. $F^vR$ is a strong
filtration, hence $G_v(R)$ is strongly graded and its degree zero
component $G_v(R)_0$ equals $\overline{R}$ which is a domain. Thus
$G_v(R)$ is a domain, hence $R$ is a domain too. Furthermore, in
this case we have $G_v(R)=\overline{R}[t,t^{-1}]$. Since
$\overline{R}$ is a Goldie domain, $\overline{R}[t,t^{-1}]$ is
also a Goldie domain and it has a skewfield of fractions
$\overline{\D}$ as well as a graded skewfield of fractions
$\overline{\D}^g$. The multiplicative set $R-\{0\}$ has
$\s(R-\{0\})=G_v(R)-\{0\}$ which is an Ore set because $G_v(R)$ is
a Goldie domain. Setting $S_0=\L-\{0\}$ we get
$\s(S_0)=G_v(R)_{\leq 0}-\{0\}$, where $G_v(R)_{\leq
0}=\dirsum_{n\leq 0}G_v(R)_n=G_v(\L)$. Clearly $\s(S_0)$ is an Ore
set in $G_v(\L)$ and $\s(S_0)^{-1}G_v(\L)=Q^g(G_v(R))$.

It follows easily that for every $p\geq 1$, $S$ maps to an Ore set
$S^{(p)}$ of $\L/m_v^p\L$. Given $s\in S_0$, $a\in R$, the left Ore condition
for $\s(S_0)$ yields an $s'\in S_0$, $a'\in R$, such that $s'a-a's\in m_v\L$,
say $s'a-a's=\pi^mb$ for some $b\in\L$, $m\in\ZZ$ (here $\pi$ is a generator
of $m_v$). Now $s''b-a''s\in m_v\L$ yields
$(s''s')a-(s''a')s=\pi^m(a''s-\pi^{m'}y$ for some $y\in\L$, $m'\neq 0$. Hence
$(s''s')a-(s''a'+\pi^ma'')s=\pi^{m+m'}y$. Consequently the microlocalization
$Q^{\mu}(\L)$ of $\L$ at $S_0$ may be defined by
$Q^{\mu}(\L)=\lim\limits_{\longleftarrow}S^{(p)^{-1}}(\L/m_v^p\L).$
Obviously $Q^{\mu}(R)=Q^{\mu}(\L)$ is a skewfield and it has a
filtration whose associated graded ring is $\overline{\D}^g$,
which is a domain and a graded skewfield. Now we can apply
Proposition \ref{ring of fractions} and conclude that the
filtration on $Q^{\mu}(R)$ is a valuation filtration. \ep

\section{Lifting from graded to filtered rings}
\noindent Let $A$ be a filtered $K$-algebra with a finite
filtration $FA$ and $\L\subset A$ a subring. Then $FA$ induces a
filtration $F\L$ on $\L$ such that $G_F(\L)\subset G_F(A)$. The
next result shows that there is a strong connection between graded
and $F$-reductors.

\bpr\label{connection} Let $A$ be as before and $\L\subset A$ a
subring.
\begin{itemize}
\item[$(i)$] If $\L\subset A$ is an $F$-reductor, then $G_F(\L)\subset
G_F(A)$ and $\widetilde{\L}\subset\widetilde{A}$ are graded
reductors.
\item[$(ii)$] If $G_F(\L)\subset G_F(A)$ or
$\widetilde{\L}\subset\widetilde{A}$ are graded reductors, then
$\L\subset A$ is an $F$-reductor.
\end{itemize} \epr
\bp (i) In order to show that $G_F(\L)\subset G_F(A)$
is a graded reductor we use Lemma \ref{lattices}(i). That
$\widetilde{\L}\subset\widetilde{A}$ is a graded reductor follows
immediately by definition.\\ (ii) By induction using Lemma
\ref{lattices}(ii). \ep

Consequently any $F$-reductor give rise to a graded reductor. On
the other side, a graded reductor is an $F$-reductor, where $FR$
is the grading filtration.

While the finite dimensional $K$-algebras have always a reductor,
see Van Oystaeyen (1975, Proposition 54 and Theorem 56), in the
infinite dimensional case the construction of reductors is not
always easy.


\bexs {\em (i) Consider $g$ a finite dimensional Lie algebra over
$K$ and $A=U(g)$ the enveloping algebra of $g$. We may define a
finite dimensional Lie algebra $g_{O_v}$ over $O_v$ with the same
basis and the induced bracket. Let us fix a $K$-basis
$\{x_1,\ldots,x_n\}$ for $g$. We have structure constants
$\l_{ij}^k\in K$ with $[x_i,x_j]=\sum_{k=1}^n\l_{ij}^kx_k$.
Without loss of generality we may assume that $\l_{ij}^k\in O_v$
(up to multiplying all $x_i$ by a suitable constant in $O_v$) but
not all in $m_v$. Set $g_{O_v}=O_vx_1+\cdots+O_vx_n$. This is a
Lie $O_v$-algebra with the induced bracket. Furthermore, $g_{O_v}$
is an $O_v$-lattice in $g$. On $\overline{g}=g_{O_v}/m_vg_{O_v}$
we define a Lie algebra structure over $k_v$ by setting
$[\overline{x}_i,\overline{x}_j]=\sum_{k=1}^n\overline{\l_{ij}^k}
\overline{x}_k$, where the $\overline{x}_i$ are the images of the
$x_i$ in $g_{O_v}$ and $\overline{\l_{ij}^k}$ are the images of
$\l_{ij}^k$ in $k_v$. By our assumptions $\overline{g}$ is not the
trivial Lie algebra. Of course, $\overline{g}$ depends on the
choice of the $K$-basis in $g$.

Let $\L=U_{O_v}(g_{O_v})$ be the enveloping algebra of $g_{O_v}$.
Consider on $A$ the standard filtration $FA$ and on $\L$ the induced
filtration $F\L$. We have that the filtration $FA$ is finite,
$G_F(\L)=O_v[X_1,\ldots,X_n]$ is a subring of the polynomial ring
$G_F(A)=K[X_1,\ldots,X_n]$ and obviously it is a graded reductor.
Proposition \ref{connection}(ii) shows that $\L$ is an $F$-reductor.
Furthermore, the filtration $F^vA$ is $\G$-separated,
$G_v(A)$ is a domain isomorphic to $U_{k_v}(\overline{g})$ and
thus we can extend $v$ to $D(g)=Q_{cl}(U(g))$. \\(ii) For the
Weyl algebra $R=A_n(K)$ we take $\L=A_n(O_v)$. We claim that
$\L$ is a graded reductor of $R$. First note that
$\L$ is a free $O_v$-module, therefore it defines a good reduction
of $R$. On the other side, we have that
$\rank_{O_v}\L_n=\dim_KR_n$ for all $n\in\NN$, since $\L$ and $R$
have the same $PBW$-basis, and so $\L$ is a graded reductor of
$R$. It is also an $F$-reductor since the associated graded ring
of $A_n(K)$ with respect to the Bernstein filtration is a polynomial
ring. }\eexs


On $G_v(A)$ we define a filtration $fG_v(A)$ by putting
$f_nG_v(A)=G_v(F_nA)$, $n\in\NN$. When $G_F(\L)\subset G_F(A)$ is
a graded reductor we have a valuation filtration $f^vG_F(A)$ on
$G_F(A)$ given by $f^v_{\g}G_F(A)=(f^v_{\g}K)G_F(\L)$, $\g\in\G$. As in the
proof of Proposition \ref{gr-reductor} we can show that the
valuation filtration $f^vG_F(A)$ on $G_F(A)$ is $\G$-separated,
provided that $F^vA$ is $\G$-separated.

\bc\label{double} Let $A$ be as before and $\L\subset A$ an
$F$-reductor such that the filtration $F^vA$ is $\G$-separated.
Then $G_v(G_F(A))$ is isomorphic to $\overline{G(A)}*\G$ and
$G_f(G_v(A))=G_v(G_F(A))$. \ec \bp We have
$G_v(G_F(A))=\dirsum_{\g\in\G}(f^v_{\g}K)G_F(\L)/(f^v_{<\g}K)G_F(\L)$
and since $G_F(\L)$ is $O_v$-flat we get
$G_v(G_F(A))=\dirsum_{\g\in\G}(f^v_{\g}K)/(f^v_{<\g}K)\otimes_{O_v}G_F(\L)$.
The latter is nothing else but
$G_v(K)\otimes_{O_v}G_F(\L)=(k_v*\G)\otimes_{O_v}G_F(\L)=\overline{G(A)}*\G$.

The fact that $G_f(G_v(A))$ equals $G_v(G_F(A))$ follows by the
general compatibility result of Hussein and Van Oystaeyen (1996,
Proposition 2.4) or can be verified directly. \ep

\bpr Let $A$ be as before and $\L\subset A$ an $F$-reductor such
that $F^vR$ is $\G$-separated. If $\overline{G(A)}$ is a domain,
then also $G_F(A)$, $G_v(A)$ and $A$ are domains. \epr \bp By
Corollary \ref{double} we get that $G_v(G_F(A))$ is a domain.
Since all principal symbol maps involved are multiplicative, we
obtain that $G_F(A)$ is a domain, hence $A$ is a domain too. From
$G_f(G_v(A))=G_v(G_F(A))$ we get that $G_f(G_v(A))$ is a domain
and thus we deduce that $G_v(A)$. \ep

Now the problem of finding an extending non-commutative valuation
is reduced to finding a reductor in an associated graded ring
having a domain for its reduction.


In the following we show that in case the algebras are given by a
finite number of relations we may relate the existence of a
reductor to properties of good reduction. Let $A$ be a filtered
$K$-algebra with a finite filtration $FA$, $n=\dim_KF_1A$ and
assume that $A=K[F_1A]$. Then we may give $A$ as an epimorphic
image of $\mathcal{F}=K<\underline{X}>$, the free $K$-algebra on
the set $\underline{X}=\{X_1,\ldots,X_n\}$. Say
$\{x_1,\ldots,x_n\}$ is a $K$-basis for $F_1A$ and
$\pi:K<\underline{X}>\longrightarrow A$ the canonical $K$-algebras
morphism given by $\pi(X_i)=x_i$, $i=1,\ldots,n$. The filtration
we use on $\mathcal{F}$ is the degree filtration (the total degree
in the $X_1,\ldots,X_n$). It follows that $\pi$ is a {\em strict}
filtered morphism, i.e. $F_n\mathcal{F}$ is taken exactly to
$F_nA$ by $\pi$. Denote by $\mathcal{R}$ the ideal of relations of
$A$, that is $\mathcal{R}=\ker\pi$. Then
\begin{equation}\label{2}
0\longrightarrow\mathcal{R}\longrightarrow\mathcal{F}\longrightarrow
A\longrightarrow 0
\end{equation}
is a strict exact sequence and therefore we obtain an exact
sequence of graded $\mathcal{F}$-modules
$0\longrightarrow\widetilde{\mathcal{R}}\longrightarrow\widetilde{\mathcal{F}}
\longrightarrow \widetilde{A}\longrightarrow 0$ by passing to the
Rees objects; see Li and Van Oystaeyen (1995). Also from strict
exactness of sequence (\ref{2}) it follows that
$G_{\mathcal{F}}(A)=G_F(R)$, where $G_{\mathcal{F}}(A)$ is the
associated graded ring of $A$ as a filtered $\mathcal{F}$-module.
Thus the sequence $0\longrightarrow G(\mathcal{R})\longrightarrow
G(\mathcal{F})\longrightarrow G(A)\longrightarrow 0$ is exact in
$G(\mathcal{F})$-gr. Since the filtration on $\mathcal{F}$ stems
from its gradation it follows that
$G(\mathcal{F})\simeq\mathcal{F}$ and under this isomorphism
$G(\mathcal{R})$ corresponds to the ideal $\mathcal{R}^{\bullet}$
in $\mathcal{F}$, where $\mathcal{R}^{\bullet}$ is the graded
ideal generated by the homogeneous components of highest degrees
of elements of $\mathcal{R}$.

Let us recall Theorem 2.13 from Hussein and Van Oystaeyen (1996).

\bt\label{lifting} Let $A$ be as before. If $G_F(A)$ reduces well
at $O_v$, then $G_F(A)$ is isomorphic to
$\mathcal{F}/\mathcal{R}^{\bullet}$ where $\mathcal{R}$ is
generated as a two-sided ideal by a finite set
$\{p_1(\underline{X}),\ldots,p_d(\underline{X})\}$ of elements of
$\mathcal{F}$ having as homogeneous parts of highest degree
$q_1(\underline{X}),\ldots,q_d(\underline{X})$ which generate
$\mathcal{R}^{\bullet}$ as a two-sided ideal. Moreover $A$ reduces
well at $O_v$. \et

Proposition \ref{connection} and Theorem \ref{lifting} together
provide complete information concerning lifting properties from
$G_F(A)$ to $A$ with respect to the existence of valuation rings
in $Q_{cl}(A)$ (or on some microlocalization in case
non-Noetherian or non-Goldie rings are being considered) extending
$O_v$. Note that for $\G=\ZZ$ Theorem \ref{lifting} allows
ramification of valuations while the reductor property involved in
Proposition \ref{connection} provides only unramified extensions
of valuations. The latter is the case in $D(g)=Q_{cl}(g)$, Weyl
algebras (see Van Oystaeyen and Willaert, 1996), Sklyanin algebras
(see Hussein and Van Oystaeyen, 1996) and many more. The
generalized gauge algebras provide an example for the first case.

\section*{Acknowledgments}
This research was supported by the bilateral project ``Hopf
Algebras in Algebra, Topology, Geometry and Physics" of the
Flemish and Romanian governments. The first author greatfully
acknowledges the support of the University of Antwerp.




\end{document}